\documentclass[12pt]{amsart}
\usepackage[T1]{fontenc}
\usepackage[utf8]{inputenc}
\usepackage{newpxtext}
\usepackage[a4paper,left=23.5mm,right=23.5mm,top=21mm,bottom=23mm,
  headheight=22pt,headsep=7mm,footskip=11mm]{geometry}
\usepackage{amsmath,amsthm}
\usepackage{newpxmath}
\usepackage{microtype,booktabs,array}
\usepackage{xcolor,tikz}
\usetikzlibrary{arrows.meta,calc}
\usepackage{caption,fancyhdr,needspace,wrapfig}
\usepackage[hidelinks]{hyperref}

\definecolor{ink}{HTML}{234D54}
\definecolor{quiet}{gray}{0.42}
\definecolor{gridgray}{gray}{0.81}
\renewcommand{\headrulewidth}{.25pt}
\renewcommand{\headrule}{\hbox to\headwidth{\color{gridgray}\leaders\hrule height\headrulewidth\hfill}}
\fancypagestyle{firstpage}{\fancyhf{}\renewcommand{\headrulewidth}{0pt}
  \fancyfoot[C]{\footnotesize\color{quiet}\rule[.45ex]{12mm}{.25pt}\enspace
  \thepage\enspace\rule[.45ex]{12mm}{.25pt}}}
\newtheoremstyle{note}{7pt}{8pt}{\itshape}{}{\bfseries\color{ink}}{.}{.6em}{}
\theoremstyle{note}
\newtheorem*{theorem}{Theorem}
\newtheorem*{proposition}{Proposition}
\newcommand{\AS}{\mathrm{AS}}
\newcommand{\DT}{\mathrm{DT}}
\newcommand{\UA}{U_{\mathrm{AS}}}
\newcommand{\UD}{U_{\mathrm{DT}}}
\newcommand{\proofstep}[1]{\par\addvspace{8pt}\needspace{3\baselineskip}
  \noindent{\bfseries\color{ink}#1.}\enspace\ignorespaces}
\tikzset{
  cellgrid/.style={draw=gridgray,line width=.25pt},
  boardedge/.style={draw=quiet,line width=.45pt},
  trace/.style={draw=ink,line width=1.05pt,line cap=round,line join=round},
  oldtrace/.style={draw=quiet,line width=.75pt,dash pattern=on 3pt off 2pt,
    line cap=round,line join=round},
  continuation/.style={draw=ink,line width=.9pt,densely dotted,line cap=round},
  celldot/.style={circle,fill=ink,inner sep=0pt,minimum size=2.3pt},
  cellname/.style={font=\scriptsize,fill=white,inner sep=1.1pt,text=black},
  panelname/.style={font=\small,text=ink},
  shadowextension/.style={draw=ink,line width=1.15pt,
    dash pattern=on .3pt off 1.6pt,line cap=round}
}
\newcommand{\fivecellgrid}{%
  \foreach \i in {0,...,5}{
    \draw[cellgrid] (\i-.5,-.5)--(\i-.5,4.5);
    \draw[cellgrid] (-.5,\i-.5)--(4.5,\i-.5);}
  \draw[boardedge] (-.5,-.5) rectangle (4.5,4.5);}

\title{Piercing a Chessboard with Restricted Slopes}
\author[Giedrius Alkauskas]{Giedrius Alkauskas}
\address{Vilnius University, Department of Mathematics and Informatics,
	Naugarduko 24, LT-03225 Vilnius, Lithuania}
\email{giedrius.alkauskas@mif.vu.lt}
\date{\today}

\begin{document}
	\begin{abstract}
Is it true that for $n\geq 3$, $n-1$ lines are needed to pierce every cell of an $n\times n$ board? This problem was independently rediscovered several times: Alkauskas (1998),  Harel-Ackerman-Pinchasi (2012, without claiming the $n-1$ bound), B\'ar\'any-Frankl (2021). However, apparently the first to pose it was T. Keller (Honolulu, Hawaii), in 1980, in the problem section of \emph{The American Mathematical Monthly}. In this paper, we give an affirmative answer to Keller's problem when all lines have slopes in $(-\infty,-1]\cup(0,1]$.
	\end{abstract}

\maketitle
\proofstep{History and the result} A line is said to \emph{pierce} a unit square if it meets its interior. Let $p_n$ be the least number of straight lines piercing all $n^2$ unit cells of an $n\times n$ board. Is it true that for $n\geq 3$, $p_n=n-1$?\\

The problem occurred to the author in 1998 (see the printed version dated 2004 \cite{MSE}, and also \cite{Alkauskas2}, Problem 9). Despite having the International Mathematical Olympiad as a final destination in mind, the author abandoned this endeavour owing to his inability to find a solution. In 2013, Ackerman and Pinchasi \cite{AckermanPinchasi2013} proved that $p_n\geq\left\lceil\frac{2}{3}n\right\rceil$. The proof has two main steps: normalization among each family and the double cover count. The authors attribute the piercing problem to Lapid Harel, an undergraduate student (2012) at the Technion. In 2023, Richter proved yet again that $p_n\geq\left\lceil\frac{2}{3}n\right\rceil$ using essentially the same argument. Kerimov extended the monotone path result to rectangular
boards \cite{Kerimov2015}. A separate development began with a question of B\'ar\'any and Frankl in 2021 \cite{BaranyFrankl2021}. Ambrus, B\'ar\'any, Frankl, and Varga proved
 $p_n>0.7n$ for sufficiently large $n$, using a linear programming technique. The authors also gave the following valuable heuristic for why $n-1$ is a plausible candidate. The symmetric plank theorem of Ball \cite{Ball1991} implies that, given any set of $n-1$ lines in the plane, there always exists a $1\times 1$ square inside the board whose sides are parallel to the grid and whose interior is not intersected by any of the lines. The discrete question is: given any set of \(n-2\) lines, can such a unit square be chosen with its bottom-left vertex at integer coordinates?\\

Due to the utmost naturalness of the problem in question, it was the deep belief of the author that such a conjecture must have been formulated well before 1998. And indeed, an extensive search revealed that this exact question was posed by T. Keller as Problem E2862 in the {\sc American Mathematical Monthly}. Most likely, an amateur mathematician Tim Keller from Hawaii.\\

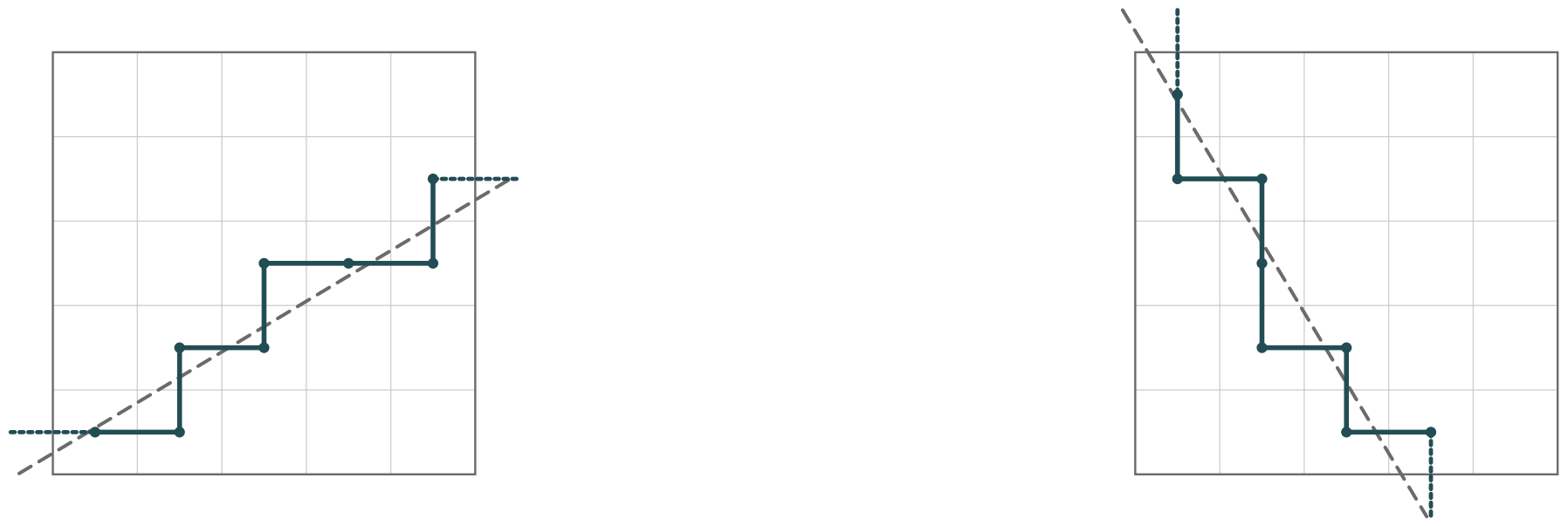
\begin{figure}[htbp]
	\centering
	\begin{tikzpicture}[x=.64cm,y=.64cm]
		\begin{scope}
			\draw[cellgrid,step=1] (0,0) grid (5,5);
			\draw[boardedge] (0,0) rectangle (5,5);
			\draw[oldtrace] (-.4,.01)--(5.4,3.49);
			\draw[continuation] (-.5,.5)--(.5,.5) (4.5,3.5)--(5.5,3.5);
			\draw[trace] (.5,.5)--(1.5,.5)--(1.5,1.5)--(2.5,1.5)--
			(2.5,2.5)--(3.5,2.5)--(4.5,2.5)--(4.5,3.5);
			\foreach \x/\y in {.5/.5,1.5/.5,1.5/1.5,2.5/1.5,2.5/2.5,
				3.5/2.5,4.5/2.5,4.5/3.5}{\node[celldot] at (\x,\y) {};}
			\node[panelname] at (2.5,-.85) {$\AS$: $0<m\le1$};
		\end{scope}
		\begin{scope}[xshift=8.2cm]
			\draw[cellgrid,step=1] (0,0) grid (5,5);
			\draw[boardedge] (0,0) rectangle (5,5);
			\draw[oldtrace] (-.15,5.5)--(3.45,-.5);
			\draw[continuation] (.5,5.5)--(.5,4.5) (3.5,.5)--(3.5,-.5);
			\draw[trace] (.5,4.5)--(.5,3.5)--(1.5,3.5)--(1.5,2.5)--
			(1.5,1.5)--(2.5,1.5)--(2.5,.5)--(3.5,.5);
			\foreach \x/\y in {.5/4.5,.5/3.5,1.5/3.5,1.5/2.5,1.5/1.5,
				2.5/1.5,2.5/.5,3.5/.5}{\node[celldot] at (\x,\y) {};}
			\node[panelname] at (2.5,-.85) {$\DT$: $m\le-1$};
		\end{scope}
	\end{tikzpicture}
	\caption{A permitted line and its path.}
	\label{fig:lines}
\end{figure}

In the current note, a partial answer to Keller's problem is given. It is the first ever occurrence of $n-1$ as an exact lower bound in the piercing problem.
\begin{theorem}
If lines with slopes in $(-\infty,-1]\cup(0,1]$ pierce every unit cell
of an $n\times n$ square, then there are at least $n-1$ of them.
\end{theorem}
\proofstep{Lines and paths}
Translate each line slightly, without losing any pierced cell, so that
it avoids grid vertices.
This is possible because there are only finitely many relevant cells
and vertices. Join the centres of the pierced cells in their order
from left to right, as in Figure~\ref{fig:lines}.
A line of slope $0<m\le1$ gives right and up steps, with no two consecutive
up steps; call this an $\AS$ path (\emph{ascending shallow}). A line of slope $m\le-1$ gives
right and down steps, with no two consecutive right steps; call this a
$\DT$ path (\emph{descending steep}). The rules also hold at entrance and exit. Outside the board, continue $\AS$ paths horizontally and $\DT$ paths vertically. A quarter-turn interchanges the two types, reversing direction when necessary. The main Theorem follows from the stronger claim:
\begin{proposition} 
Any covering of the $n\times n$ board with $\mathrm{AS}$-$\mathrm{DT}$ paths has at least $n-1$ members.
	\end{proposition}
\begin{wrapfigure}{l}{0.58\textwidth}
\vspace{-6pt}
	\centering
	\includegraphics[width=0.97\linewidth]{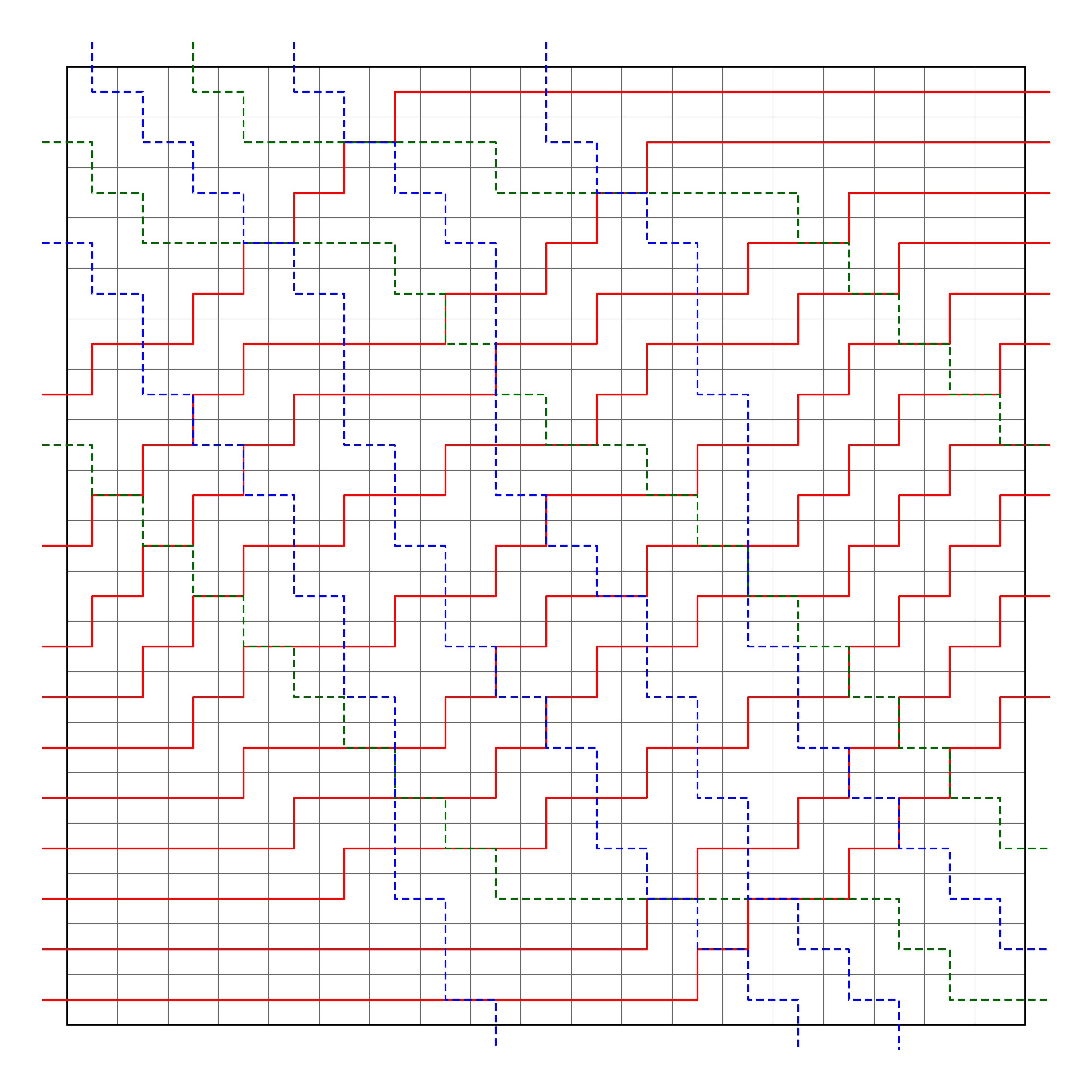}
	\caption{$n=19$ is the smallest integer such that an $n\times n$ board can be covered with $n-2$ paths of 3 types:  10AS+3DS+4DT.}
	\label{3t}
\vspace{-7pt}
\end{wrapfigure}

\noindent By analogy, we define the remaining two path classes: AT (right and up steps, no two right in a row) and DS paths (down and right, no two down in a row). In essence, this abstraction is a refinement of the Ackerman and Pinchasi setting of monotone paths. The separation property below is a direct analogue of their Lemma 1. However, it is still impossible to settle Keller's problem in full generality with such a  refinement. For, even if we restrict ourselves to three-type covers, a $19\times 19$ board can already be covered with $17$ paths (Figure \ref{3t}). More generally, a linear drop manifests itself as $n\rightarrow\infty$: an $n\times n$ board can be covered by no more than $0.983 n$ paths of three types for $n\gg 0$.
 
 In his original solution, T. Keller proposed one particular $(n-1)$-piercing of an $n\times n$ board. Ambrus-B\'ar\'any-Frankl-Varga gave a construction of a completely different family. These examples involve exactly two line classes and therefore the theorem applies. However, the variety of $(n-1)$-piercings of an $n\times n$ board is much bigger, and it includes piercings using $3$ and all $4$ classes. Thus, in our opinion, there is still not enough evidence that $n-1$ is the correct lower bound, and the possibility of a counterexample for $n\gg 0$ is not excluded. \footnote{The exact constant for $4$ type path cover is currently under investigation.}

\proofstep{Intersection and separation}
 Identify the board with $\{0,\ldots,n-1\}^2$, where $(0,0)$ is the bottom-left cell.
An $\AS$ path and a $\DT$ path have at most two common
cells. Their common cells lie in one row or one column.

\begin{figure}[htbp]
	\centering

\begin{tikzpicture}[x=5mm,y=5mm,
	baseline=(current bounding box.center)]
	
	\newcommand{\drawboard}{
		\draw[black, semithick] (0,0) rectangle (8,8);
		\foreach \x in {1,...,7}
		\draw[black, semithick] (\x,0)--(\x,8);
		\foreach \y in {1,...,7}
		\draw[black, semithick] (0,\y)--(8,\y);
	}
	
	\newcommand{\singledot}[2]{%
		\fill[black!70] ({#1+0.5},{#2+0.5})
		circle (2.2pt);
	}
	
	\newcommand{\doubledot}[2]{%
		\fill[black!70] ({#1+0.37},{#2+0.63})
		circle (2.0pt);
		\fill[black!70] ({#1+0.63},{#2+0.37})
		circle (2.0pt);
	}
	
	\begin{scope}
		\drawboard
		
		\singledot{0}{7}
		\doubledot{1}{6}
		\singledot{2}{5}
		\doubledot{5}{2}
		\singledot{6}{1}
		
		\node at (4,-0.6) {$\alpha$};
	\end{scope}
	
	\node at (9.6,4) {$\Longrightarrow$};
	
	\begin{scope}[shift={(11.2,0)}]
		\drawboard
		
		\singledot{0}{7}
		\singledot{1}{6}
		\singledot{2}{5}
		\singledot{3}{4}
		\singledot{5}{2}
		\singledot{6}{1}
		\singledot{7}{0}
		
		\node at (4,-0.6) {$b$};
	\end{scope}
	
\end{tikzpicture}

\caption{The stretch from $\alpha$ to $b$}
\label{norml}
\end{figure}
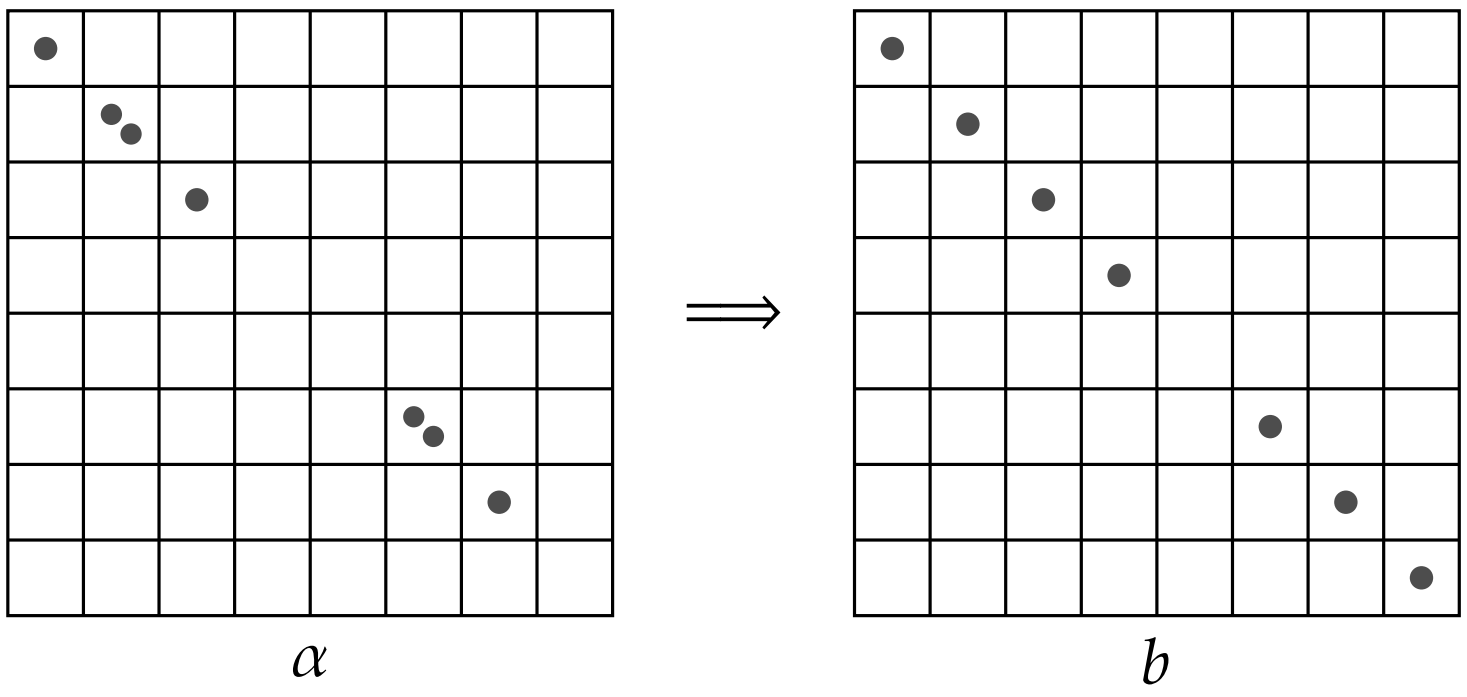

A family of $m$ paths of one type can be replaced by $\leq m$ pairwise disjoint,
globally ordered paths of that type, without losing any covered
cell. On
$x+y=t$, the $i$th full path has a single cell $\left(a_i(t),t-a_i(t)\right)$, where
\begin{eqnarray}
  0\le a_i(t+1)-a_i(t)\le1,
  \quad a_i(t+2)\ge a_i(t)+1.
  \label{kelias}
\end{eqnarray}
For all $t$, sort
$\{a_1(t),\ldots,a_m(t)\}$
as $\alpha_1(t)\le\cdots\le\alpha_m(t)$. Replacing the $a_i$ by the $\alpha_i$ preserves \eqref{kelias}. Indeed, if two AS paths have two meetings between which their relative order is reversed, the intervening fragments can be interchanged. It is easy to check that the AS rules are preserved. Further, let $b_i(t)=\max_{j\le i}\{\alpha_j(t)+i-j\}$. Equivalently, on each diagonal we process the points
from left to right, placing each at the first available
position at or to the right of its original position.
Every originally occupied position remains occupied.  Figure~\ref{norml} illustrates this in a generic case.
Since pointwise maxima preserve \eqref{kelias},
the $b_i$ are valid AS paths. Moreover,
$b_{i+1}(t)\ge b_i(t)+1$, so these paths are
disjoint and globally ordered. Restrict them to
the board and discard those with empty traces (the \emph{trace} of a path is its intersection with the board). The DT case follows by symmetry.

\proofstep{An extremal covering}
Suppose the theorem fails. Choose the smallest $n$ admitting a
covering by at most $n-2$ paths, $p$ of them
being AS and $q$ being DT. Obviously, $n\geq 4$.
Among these coverings choose one maximizing $|\UA|+|\UD|$, where
$\UA,\UD$ are the two covered-cell unions. Separation cannot decrease
this sum, so we may assume that each family is pairwise disjoint and ordered.

No path has its entire trace contained in the union of two adjacent sides (columns $0,n-1$ and rows $0,n-1$). Indeed, deleting that path together with the
corresponding sides would give a cover of an $(n-1)\times (n-1)$ square with
at most $n-3$ paths, contrary to minimality.
A $\DT$ path through the corner $(0,0)$ would lie in the union of the west and south sides; analogously for a DT path passing through $(n-1,n-1)$.
Consequently, $(0,0)$ and $(n-1,n-1)$ are covered by the $\AS$ family only, and $(0,n-1)$, $(n-1,0)$ are covered by 
 the $\DT$ family only. In particular, $p,q\geq 1$.

\begin{figure}[tbp]
	\centering
	\begin{tikzpicture}[x=.64cm,y=.64cm]
		\foreach \s/\ttl in {0/before,8.7/after}{
			\begin{scope}[shift={(\s,0)}]
				\foreach \i in {0,...,6}{\draw[cellgrid] (\i-.5,-.75)--(\i-.5,3.5);}
				\foreach \j in {0,...,4}{\draw[cellgrid] (-.5,\j-.5)--(5.5,\j-.5);}
				\draw[boardedge] (-.5,3.5)--(5.5,3.5);
				\node[panelname] at (2.5,4.2) {\textsc{\ttl}};
		\end{scope}}
		\draw[oldtrace,-{Stealth[length=3pt]}]
		(1,3.92)--(1,3)--(2,3)--(2,2)--(3,2)--(3,1)--(4,1)--(4,0)--(4,-.7);
		\foreach \j in {0,...,3}{
			\node[cellname] at ({1+\j},{3-\j}) {$U_{\j}$};
			\node[cellname,text=ink] at (\j,{3-\j}) {$V_{\j}$};}
		\foreach \j in {0,...,2}{\node[cellname] at ({2+\j},{3-\j}) {$T_{\j}$};}
		\begin{scope}[shift={(8.7,0)}]
			\draw[trace,-{Stealth[length=3pt]}]
			(0,3.92)--(0,3)--(1,3)--(1,2)--(2,2)--(2,1)--(3,1)--(3,0)--(4,0)--(4,-.7);
			\foreach \j in {0,...,3}{
				\node[cellname] at ({1+\j},{3-\j}) {$U_{\j}$};
				\node[cellname,text=ink] at (\j,{3-\j}) {$V_{\j}$};}
			\foreach \j in {0,...,2}{\node[cellname,text=quiet] at ({2+\j},{3-\j}) {$T_{\j}$};}
		\end{scope}
	\end{tikzpicture}
	\caption{The north-side exchange, with $V_0=B$ and $U_0=G$.
		The $e$ lost cells $T_j$ are $\AS$-covered;
		the $e+1$ gained cells $V_j$ were $\DT$-free.}
	\label{fig:exchange}
\end{figure}
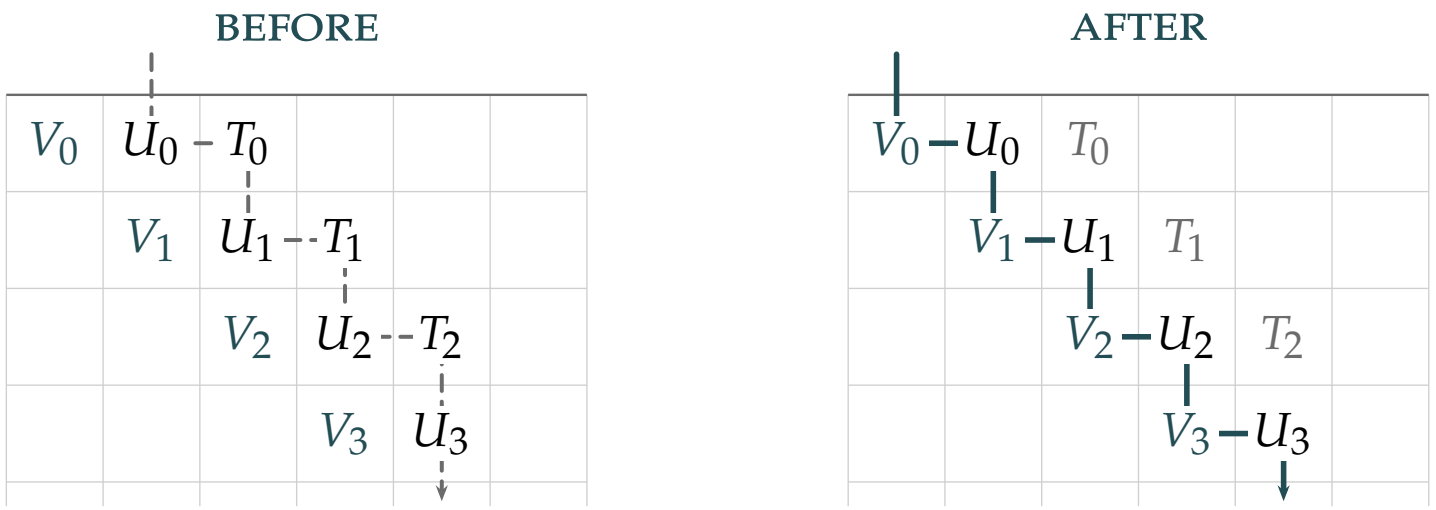

\proofstep{Ordering the boundary} We claim that the $\DT$-covered cells on the north side form a single interval starting from cell $(0,n-1)$. Assume the contrary. Let two cells satisfy
$B=(h,n-1)\notin\UD$ and $G=(h+1,n-1)\in\UD$; see Figure~\ref{fig:exchange}.
The trace of the $\DT$ path $\mathscr{P}$
through $G$ starts at $G$, since its left
neighbour $B$ is DT-free. Write the path's maximal initial alternating part as
$U_0,T_0,U_1,T_1,\ldots,T_{e-1},U_e$, with
\[
  U_j=(h+1+j,n-1-j),\quad
  T_j=(h+2+j,n-1-j),\quad
  V_j=(h+j,n-1-j).
\]
Here $V_0=B$ and $U_0=G$.
This alternating part cannot reach the bottom row, since $h+1\geq 1$. Thus after $U_e$, $\mathscr{P}$ goes down or exits east. Every $V_j$ is DT-free. This is clear for $j=0$.
For $j\ge1$, a DT path entering $V_j$ from above
would use $U_{j-1}$, already occupied by
$\mathscr{P}$. If it entered $V_j$ from the left, DT rules would force
it to pass through $V_{j-1}$, which is DT-free
by induction. Thus all $V_j$ are AS-covered.
They have distinct AS owners.

We next show that all cells $T_j$, $0\le j<e$,
are AS-covered. If $e=0$, there is nothing to prove.
For $1\le j<e$, the AS owner of $V_j$, two steps
later, must reach $T_{j-1}$ or $T_j$, both of which
are inside the board. If the AS path through \(V_0\) does not reach \(T_0\), modify it by continuing it horizontally along the north side after \(V_0\). This preserves
all its previously covered board cells and makes
it cover $T_0$. Apply the separation lemma to the modified AS
family. Every previously covered cell remains
covered, and the resulting paths are disjoint.
After separation, the AS owner of \(V_0\) reaches \(T_0\). Hence the AS owner of \(V_1\) must reach \(T_1\), rather than \(T_0\), two steps later. Inductively, the owners of
$V_2,\ldots,V_{e-1}$ must reach
$T_2,\ldots,T_{e-1}$, respectively. Thus all the cells $T_0,\ldots,T_{e-1}$ are
AS-covered in the modified covering. If any
of them had originally been AS-free, the value
$|\UA|+|\UD|$ would have increased, contrary
to maximality. Hence they were all AS-covered
in the original covering.

Now replace the initial part of the path $\mathscr{P}$ by
$V_0,U_0,V_1,U_1,\ldots,V_e,U_e$, entering $B$ from above.
The new word alternates right and down. Retain the old tail after $U_e$ if its next
step is down; if it exits east, continue the new path vertically down
column $n-1$. In counting $|\UD|$, we lose the $e$ cells $T_j$,
already $\AS$-covered, and gain the $e+1$ previously $\DT$-free cells
$V_j$. Full coverage is nevertheless preserved.
Separating the $\DT$ family preserves this increase, a contradiction.

The claim about the initial interval in row $n-1$ therefore follows. Rotation gives the exact boundary intervals:
column $x=0$ is $\AS$-covered precisely for $y\le u$, column $x=n-1$ for $y\ge z$;
row $y=n-1$ is $\DT$-covered precisely for $x\le c$, row $y=0$ for $x\ge b$.
The corner information gives $0\le u,c\leq n-2$ and $1\leq z,b\le n-1$.

\proofstep{Four shadows and extended shadows}
Let $A_{NW},A_{SE}$ be the first and last $\AS$ traces in increasing
horizontal order on $x+y=t$. The first meets the west and north sides:
its horizontal coordinate is at most $0$ at $t=0$, and at most $n-1$
at $t=2n-2$, since the southwest and northeast corners are $\AS$-covered. The last cell on the west side covered by $A_{NW}$ is exactly $(0,u)$.
The trace $A_{SE}$ similarly meets the south and east sides; its first east-side contact is $(n-1,z)$. Rotation gives the extreme $\DT$ tracks $D_{SW},D_{NE}$ meeting
the west and south sides, and the north and east sides, respectively (if $p=1$ or $q=1$, the extremes of a one-path family coincide). Moreover, $(b,0)$ is the first south contact of $D_{SW}$, and $(c,n-1)$ is the last north contact of $D_{NE}$.

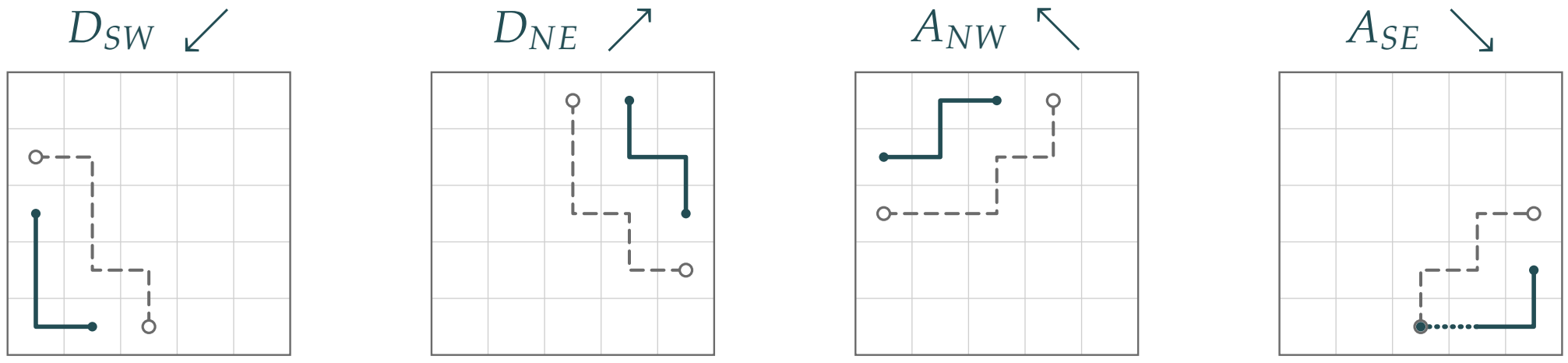
\begin{figure}[htbp]
\centering
\begin{tikzpicture}[x=.48cm,y=.48cm]
  \begin{scope}
    \fivecellgrid
    \draw[oldtrace] (0,3)--(1,3)--(1,2)--(1,1)--(2,1)--(2,0);
    \draw[trace] (0,2)--(0,1)--(0,0)--(1,0);
    \foreach \x/\y in {0/3,2/0}{\draw[draw=quiet,fill=white,line width=.6pt] (\x,\y) circle[radius=1.5pt];}
    \foreach \x/\y in {0/2,1/0}{\node[celldot] at (\x,\y) {};}
    \node[panelname] at (2,5.25) {$D_{SW}$\; $\swarrow$};
  \end{scope}
  \begin{scope}[shift={(7.5,0)}]
    \fivecellgrid
    \draw[oldtrace] (2,4)--(2,3)--(2,2)--(3,2)--(3,1)--(4,1);
    \draw[trace] (3,4)--(3,3)--(4,3)--(4,2);
    \foreach \x/\y in {2/4,4/1}{\draw[draw=quiet,fill=white,line width=.6pt] (\x,\y) circle[radius=1.5pt];}
    \foreach \x/\y in {3/4,4/2}{\node[celldot] at (\x,\y) {};}
    \node[panelname] at (2,5.25) {$D_{NE}$\; $\nearrow$};
  \end{scope}
  \begin{scope}[shift={(15,0)}]
    \fivecellgrid
    \draw[oldtrace] (0,2)--(1,2)--(2,2)--(2,3)--(3,3)--(3,4);
    \draw[trace] (0,3)--(1,3)--(1,4)--(2,4);
    \foreach \x/\y in {0/2,3/4}{\draw[draw=quiet,fill=white,line width=.6pt] (\x,\y) circle[radius=1.5pt];}
    \foreach \x/\y in {0/3,2/4}{\node[celldot] at (\x,\y) {};}
    \node[panelname] at (2,5.25) {$A_{NW}$\; $\nwarrow$};
  \end{scope}
  \begin{scope}[shift={(22.5,0)}]
    \fivecellgrid
    \draw[oldtrace] (2,0)--(2,1)--(3,1)--(3,2)--(4,2);
    \draw[trace] (3,0)--(4,0)--(4,1);
    \draw[shadowextension] (2,0)--(3,0);
    \foreach \x/\y in {2/0,4/2}{\draw[draw=quiet,fill=white,line width=.6pt] (\x,\y) circle[radius=1.5pt];}
    \foreach \x/\y in {2/0,4/1}{\node[celldot] at (\x,\y) {};}
    \node[panelname] at (2,5.25) {$A_{SE}$\; $\searrow$};
  \end{scope}
\end{tikzpicture}
\caption{The four shadows and their extensions.}
\label{fig:shifts}
\end{figure}

For each of these four traces, construct its \emph{shadow}. That is, take the portion of the path from its last contact with the first named
side to its first contact with the second, and translate it outward as in Figure~\ref{fig:shifts}. Upon restriction to the board, the translated portion loses exactly its two endpoints.
The original portion avoids the intervening corner and therefore has at least
three cells, so the shadow is nonempty. Each cell of every shadow is covered solely by paths of the opposite family.

Now extend the southwest $\DT$ shadow upward along column $x=0$ to $(0,u)$ and the
northeast $\DT$ shadow downward along column $x=n-1$ to $(n-1,z)$.
Extend the northwest $\AS$ shadow rightward along the north side to $(c,n-1)$ and the
southeast $\AS$ shadow leftward along the south side to $(b,0)$.
Call the resulting segments the \emph{extended shadows}.
The boundary intervals established above ensure that every cell added to a
shadow is also covered by the opposite family. The extensions also preserve the forbidden-step
conditions at the joins. 

For example, the last step of $D_{SW}$ is $(b,1)\to(b,0)$, so its
shadow ends at $(b-1,0)$. Its extended shadow starts at
$(0,u)$ and has $u+b$ cells. Rotation gives the other three cases.
Each $\DT$ extended shadow is $\AS$-covered and meets each of the $p$ $\AS$
paths in at most two cells, so it has at most $2p$ cells.
Similarly, each $\AS$ extended shadow has at most $2q$ cells.

On the other hand, the four extended shadows, viewed as chains of cells, form a closed walk whose successive portions run northeast, southeast, southwest, and northwest. Its total horizontal and vertical travel are each \(2(n-1)\), so its length is \(4(n-1)\). Therefore, \(4(n-1)\le 4p+4q\), the desired contradiction.

\end{document}